\documentclass{article}
\usepackage{graphicx}

\newtheorem{theorem}{Theorem}

\newtheorem{prop}[theorem]{Proposition}

\newtheorem{definition}[theorem]{Definition}

\title{Deforming Diamond}
\author{Ciprian S. Borcea and Ileana Streinu}
\date{}

\begin{document}
\maketitle

\begin{abstract}
For materials science, diamond crystals are almost unrivaled for hardness and a range of other properties. Yet, when simply abstracting the carbon bonding structure as a geometric bar-and-joint periodic framework, it is far from rigid. We study the geometric deformations of this type 
of framework in arbitrary dimension $d$, with particular regard to the volume variation of a unit cell.  
  
\end{abstract}



\section*{Introduction}
   
In this paper we survey a number of topics in the  deformation theory of periodic
bar-and-joint frameworks by investigating diamond frameworks, that is, generalizations to arbitrary dimension $d$ of the ideal atom-and-bond structure of diamond crystals. 

\medskip \noindent
Our setting is purely geometric and does not involve any physical assumptions or properties. 
A general deformation theory of periodic frameworks was introduced in \cite{BS1} and further developed in \cite{BS2, BS3}. 
Although crystallography, solid state physics and materials science have studied
crystal structures for a long time, abstract mathematical formulations are of relatively recent date \cite{Su2}. Of course, there are historical roots in problems related to sphere packings \cite{K,G,V},   lattice theory \cite{Br, M} and crystallographic groups \cite{F,Sch,Bb}.

\begin{figure}[h]
 \centering
 {\includegraphics[width=0.85\textwidth]{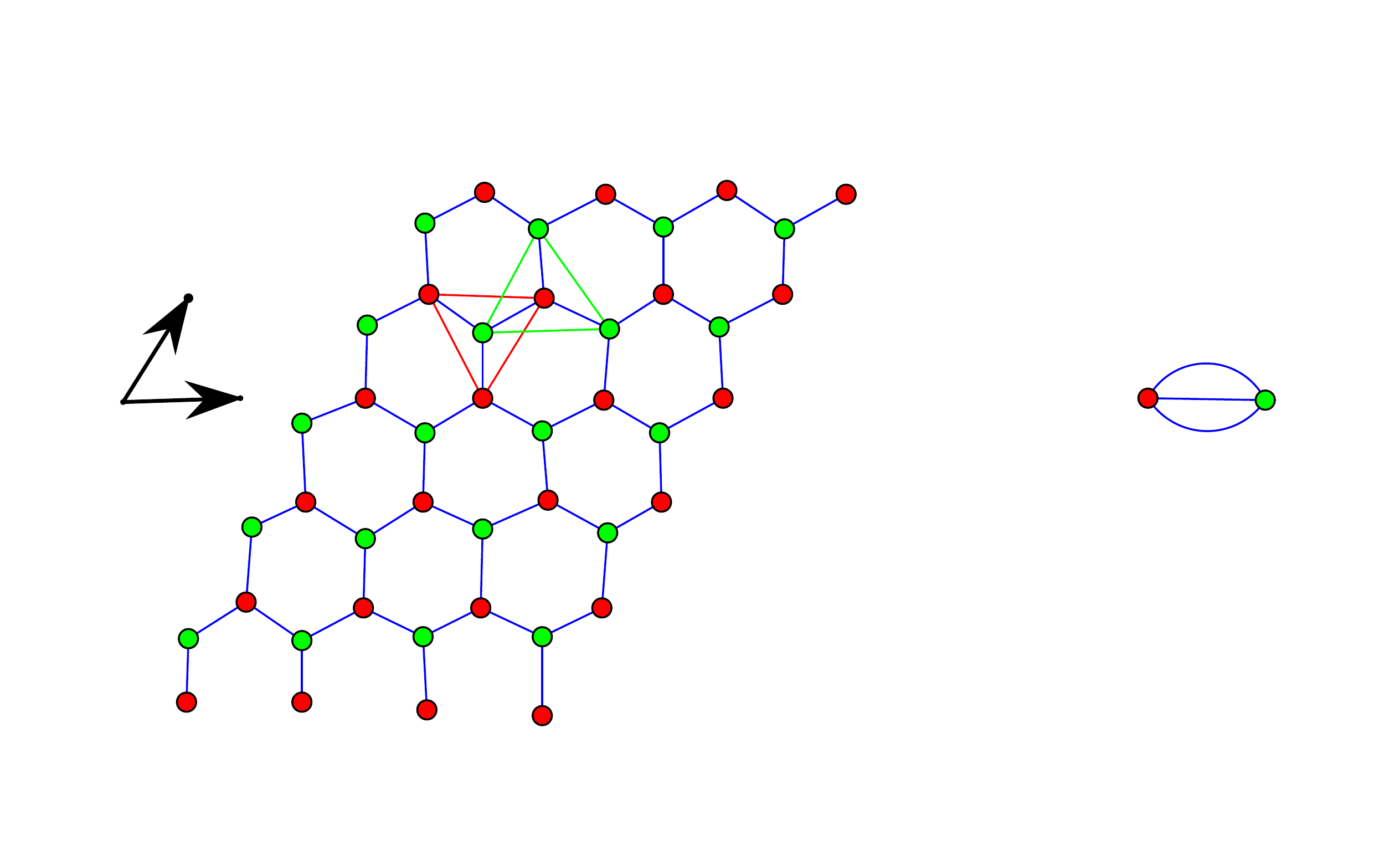}}
 \caption{ The standard diamond framework in dimension two. The generators of the periodicity lattice are emphasized by arrows. The quotient graph has two vertices connected by three edges.}
 \label{Fig2diamond}
\end{figure}

\begin{figure}[h]
 \centering
 {\includegraphics[width=0.60\textwidth]{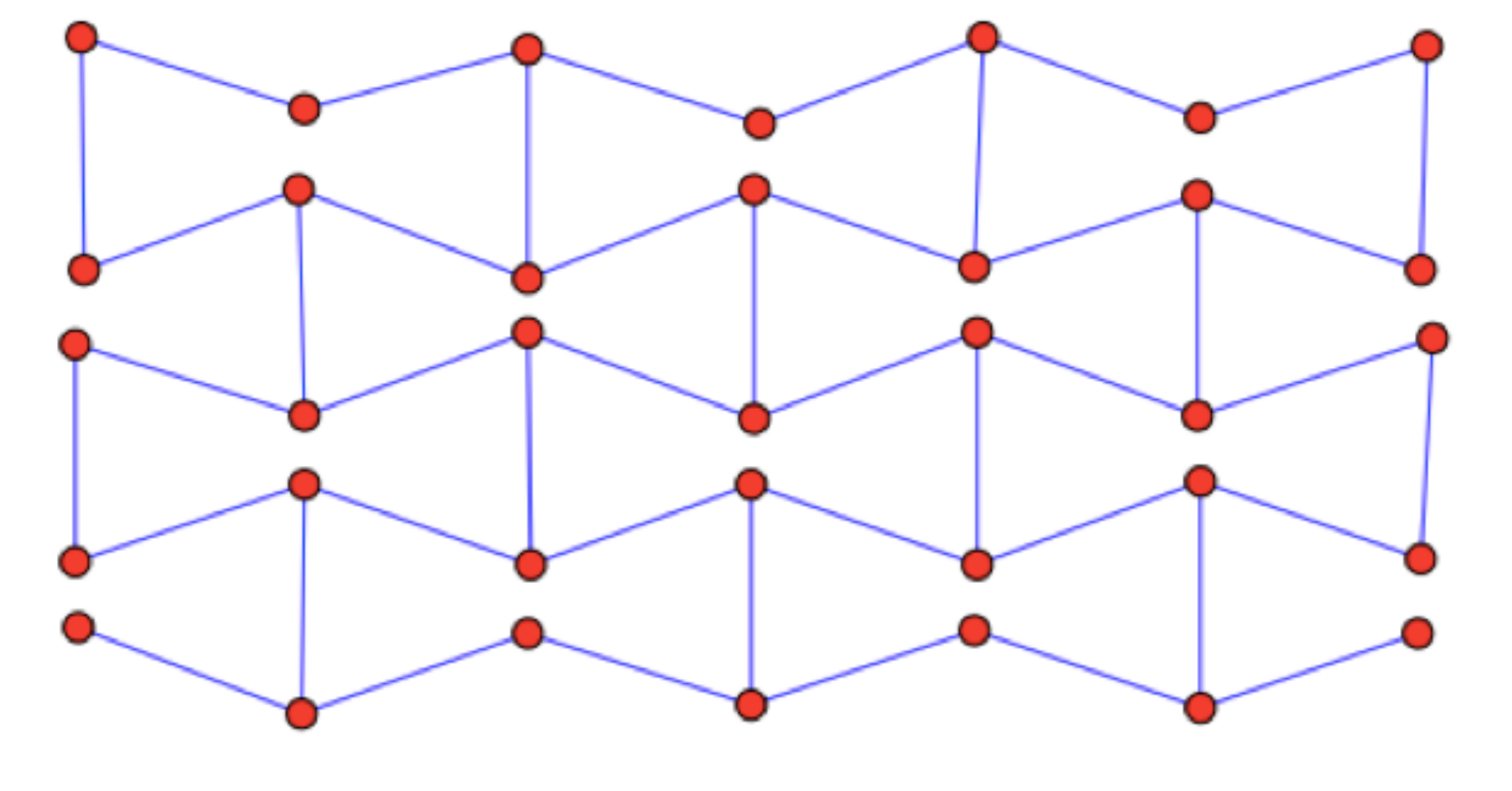}}
 \caption{ A reentrant honeycomb is a deformed diamond framework with auxetic capabilities.  A horizontal streching detrmines a vertical expansion.}
 \label{Figauxetic2D}
\end{figure}

\medskip \noindent
Diamantine frameworks have elementary geometrical descriptions in any dimension $d$. The standard planar case $d=2$ is the familiar {\em hexagonal honeycomb} illustrated in Figure~\ref{Fig2diamond}.  Another configuration, the so-called {\em reentrant honeycomb} illustrated in Figure~\ref{Figauxetic2D}, was recognized in materials science as a structure with 
negative Poisson's ratio \cite{Kol, L} and has become an emblem of auxetic behaviour \cite{ENHR, GGLR}. It will be seen below that auxetic deformation paths can be defined in strictly
geometric terms \cite{BS5} and resemble causal-lines in special relativity. 
 
\medskip \noindent
Our treatment is mostly self-contained and may serve as an ``introduction by example"  to several topics of general interest: topology of the deformation space, variation of volume per unit cell and critical configurations, possibilities for auxetic trajectories.

\section{The standard diamond framework in $R^d$}\label{sDF}

\begin{figure}[h]
 \centering
 {\includegraphics[width=0.60\textwidth]{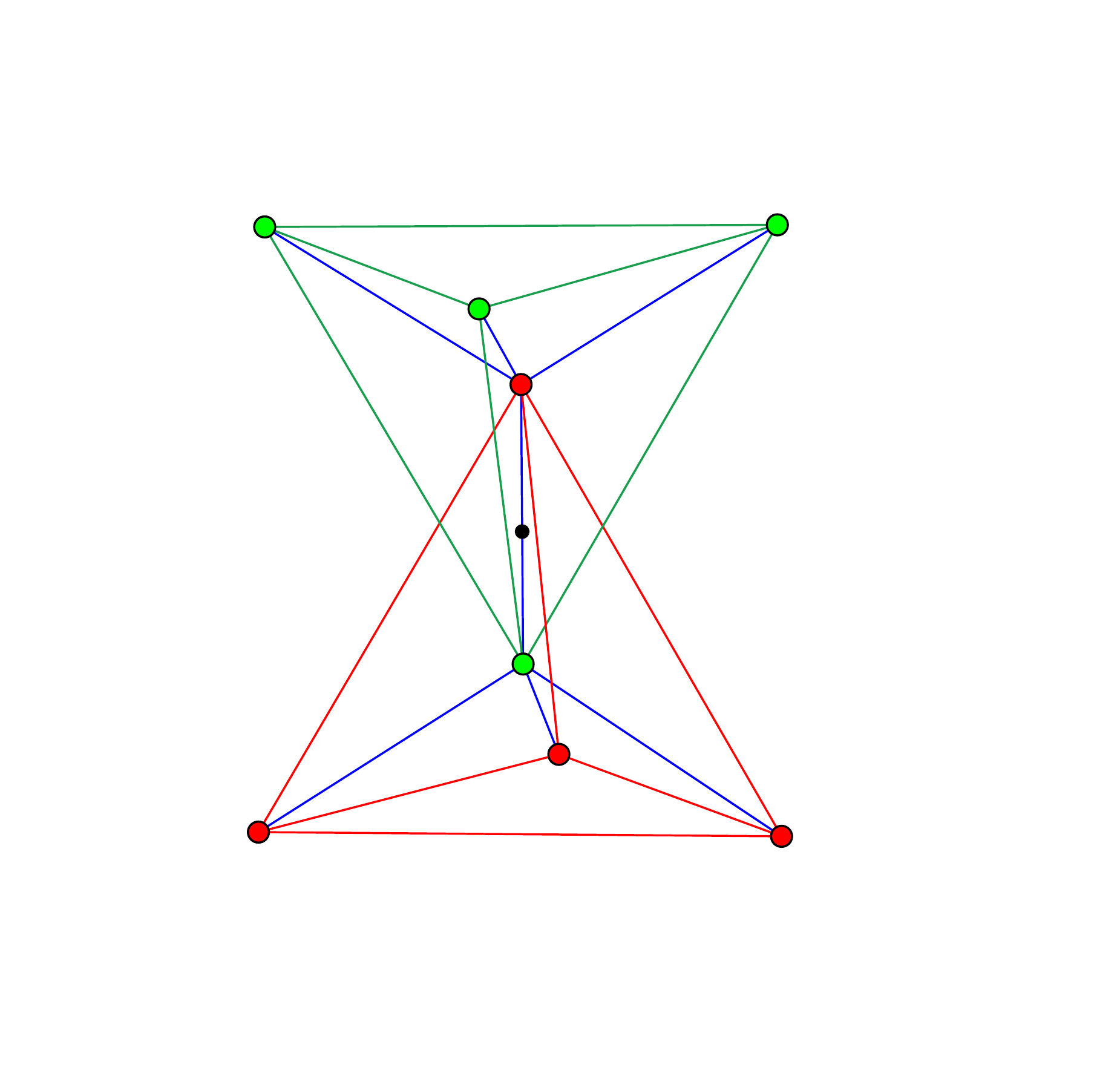}}
 \caption{ An illustration for $d=3$. The standard diamond framework uses the vertices of two
regular tetrahedra. Each tetrahedron center is a vertex of the other tetrahedron and the pair is
cetrallly symmetric with respect to the midpoint of the two centers. Centers are seen as joints connected by bars to the vertices of the other tetrahedron. By periodicity under the 
lattice of translations generated by all edge-vectors of the tetrahedra,  one obtains the full
bar-and-joint diamond framework.}
 \label{Fig3diamond}
\end{figure}

The standard or canonical diamond framework in $R^d$ can be described by starting with
a regular simplex $P_0P_1...P_d$ centered at the origin $O$. Then, we take the midpoint $M_0$
of the segment $OP_0$ and denote by $Q_0Q_1...Q_d$ the simplex obtained from the original one by central symmetry with center $M_0$. Figure~\ref{Fig3diamond} shows this setting in dimension three.
By using as {\em periodicity lattice} $\Lambda$ the rank $d$
lattice generated by the vectors $\lambda_i=P_i-P_0$, $i=1,...,d$, we can produce a periodic set of vertices which includes the vertices $P_i$ and $Q_j$, by considering all translates of $P_0$ and $Q_0=O$. For edges we take all segments $OP_i=Q_0P_i$, $i=0,...,d$ and their translates under $\Lambda$. Up to isometry and rescaling, the resulting framework is unique.

The abstract infinite graph together with the marked automorphism group represented by $\Lambda$ gives a $d$-periodic graph as defined in \cite{BS1, BS2}. There are two orbits
of vertices and $d+1$ orbits of edges under this periodicity group. In other words, the quotient
multigraph consists of two vertices connected by $d+1$ edges.

\medskip
The same $d$-periodic graph allows other placements. One may start with an arbitrary simplex
$P_0P_1...P_d$ with the origin $O$ in the interior. Then $O$ is connected with the vertices $P_i$
and provides $d+1$ edge representatives. By using the periodicity lattice $\Lambda$ generated
as above by the edge vectors of the initial simplex, one extends the finite system of vertices and
edges to a full $\Lambda$-periodic system. The corresponding frameworks will be called {\em diamantine frameworks}, frameworks of diamond type or simply diamond frameworks.

\medskip
The three dimensional standard  framework corresponds with the bonded structure of carbon atoms in (ideal) diamond crystals. The two dimensional version may be imagined as a graphene layer.

\medskip
What distinguishes the standard diamond framework among its diamantine relatives is its
full symmetry group, which is maximal. We refer to \cite{KS, Su1, Su2} for the corresponding theory of canonical placements of periodic graphs. Related aspects can be found in \cite{DF, E, BS3}. Our present inquiry will assume that a diamantine framework in $R^d$ has been given and
will be concerned with its {\em deformations} as a periodic bar-and-joint structure. Deformations will have new placements for the vertices (i.e. joints) of the structure, but in such a way that all edges preserve their length (i.e. behave like rigid bars) and the deformed framework
remains periodic with respect to the abstract periodicity group marked at the outset \cite{BS1}.
It is important to retain the fact that the representation of this abstract periodicity group by
a lattice of translations of rank $d$ is allowed to vary in deformations.

\section{Deformations}

Let $p_0,...p_d$ be the {\em edge vectors} emanating from the origin. Their squared norm will be denoted by $\langle p_i,p_i\rangle =s_i$. The lattice of periods is generated by
$p_i-p_0,\  i=1,...,d$ and with this ordering, the oriented volume of the fundamental parallelepiped (unit cell) is given by 

\begin{equation}\label{volume}
 V_I(p)=det \left[\begin{array}{llcl} 1 & 1 & ... & 1 \\ p_{i_0} & p_{i_1} & ... & p_{i_d} \end{array}\right] =det
\left[\begin{array}{ccc} p_{i_1}-p_{i_0} & ... & p_{i_d}-p_{i_0} \end{array} \right]
\end{equation}

\medskip \noindent
The deformation space of a diamond framework defined by $p$, with $V(p)\neq 0$ can be parametrized by the connected component of $p$ in   

$$ \prod_{i=0}^d S^{d-1}_{\sqrt s_i}  \setminus \{ p : V(p)=0 \}   $$

\noindent
modulo the natural action of $SO(d)$. As usual, $S^{d-1}_r$ denotes the $(d-1)$-dimensional sphere of radius $r$.

\medskip \noindent
We investigate first the {\em critical points} of $V(p)$ on the indicated product of spheres.
At a critical point $p$, we must have 

$$ \frac{\partial V}{\partial p_i}(p)=\lambda_i p_i $$

\noindent
for some Lagrange multipliers $\lambda_i$, $i=0,...,d$. By inner product with $p_j$, we obtain:

\begin{equation}
\langle \frac{\partial V}{\partial p_i}(p) , p_j \rangle = \lambda_i \langle p_i,p_j \rangle
\end{equation}

\noindent
and for $i\neq j$ this gives:

\begin{equation}\label{distinct}
\lambda_i \langle p_i,p_j \rangle = (-1)^{i+1} det| p_0 ...\hat{p}_i... p_d | 
\end{equation}
 
\noindent
where $\hat{p}_i$ means the absence of that column.

\medskip \noindent
For $i=j$ we find:

\begin{equation}\label{same}
\lambda_i|p_i|^2=V(p)+(-1)^{i+1} det | p_0 ... \hat{p}_i ... p_d |
\end{equation}

\medskip \noindent
By summation or by Euler's formula for the homogeneous function $V(p)$, we obtain:

\begin{equation}
\sum_{i=0}^d \lambda_i |p_i|^2=dV(p)
\end{equation}

\medskip \noindent
Equations (\ref{distinct}) and (\ref{same}) show that for non-zero critical values, all $\lambda_i\neq 0$ and then it follows that all inner products $\langle p_i,p_j\rangle$,
$i\neq j$ must be equal. If we denote by $\alpha$ their common value, we see from the
Gram determinant of the vectors $p_i$ that $\alpha$ must be a root of the degree  $(d+1)$
equation

\begin{equation}\label{Gram}
det \left( \begin{array}{lllll} s_0 & \alpha & \alpha & ... & \alpha \\
                                             \alpha & s_1 & \alpha & ... & \alpha \\
                                              .. & .. & .. & ... & .. \\
                                              \alpha & \alpha & \alpha & ... & s_d \end{array} \right) =0
\end{equation} 

\medskip \noindent
For $\alpha=s_i=|p_i|^2$, the determinat on the left hand side of (\ref{Gram}) equals
$s_i \prod_{j\neq i} (s_j-s_i)$. If we assume distinct values $s_i$ and the natural ordering
$s_0 < s_1 < ... < s_d$, we see that equation (\ref{Gram}) must have a root in every interval
$(s_i,s_{i+1})$, since the determinant has opposite signs at the endpoints. Thus, in the general case of distinct $s_i$, equation (\ref{Gram}) has $d$ positive simple roots. The remaining root
must be negative since the term of degree one in $\alpha$ is zero i.e. the sum of all products of $d$ roots vanishes. By continuity,  we still have {\em a unique negative root} $\alpha_{-}$ when some of the values $s_i$ become equal.

\medskip \noindent
With respect to {\em positive} roots of  (\ref{Gram}),  we note that the inequality $\langle p_i,p_j\rangle \leq |p_i|\cdot |p_j|=(s_is_j)^{1/2}$,
implies that only the {\em smallest among them} $\alpha_{+}$ may correspond to a cricical point. Moreover,
from $V(p)\neq 0$, we infer that this root must be simple. Indeed, otherwise, with t$s_0=s_1\leq s_2\leq ...\leq s_d$, the root must be $s_0=s_1=\langle p_0,p_1\rangle$ and this forces $p_0=p_1$,  resulting in $V(p)=0$.  In fact, for the simple root $\alpha_{+}$ one may confirm
the  more precise location $\alpha_{+}\in (s_0, s_0^{1/2}s_1^{1/2})$ by evaluating the
determinant in (\ref{Gram}) at $\alpha=s_0^{1/2}s_1^{1/2}$ and finding it negative.
With $r_i=s_i^{1/2}$ and $\alpha=r_0r_1$, we have

 $$ det \left( \begin{array}{lllll} s_0 & \alpha & \alpha & ... & \alpha \\
                                             \alpha & s_1 & \alpha & ... & \alpha \\
                                              .. & .. & .. & ... & .. \\
                                              \alpha & \alpha & \alpha & ... & s_d \end{array} \right)  =
s_0s_1det \left( \begin{array}{lllll} 1 & 1 & r_1 & ... & r_1\\
                                            1 & 1 & r_0 & ... & r_0 \\
                                             r_1 & r_0 & s_2 & \alpha ... \\
                                              .. & .. & .. & ... & .. \\
                                              r_1 & r_0 & \alpha & ... & s_d \end{array} \right) =  $$

$$ = s_0s_1(r_0-r_1)^2 det \left( \begin{array}{lllll} 0 & 1 & 1 & ... & 1 \\
                                             1 & s_2 & \alpha & ... & \alpha \\
                                              .. & .. & .. & ... & .. \\
                                             1 & \alpha & \alpha & ... & s_d \end{array} \right) = $$

$$ = -s_0s_1(r_0-r_1)^2 \prod_{i=2}^d (s_i-\alpha)(\sum_{i=2}^d \frac{1}{s_i-\alpha}) < 0. $$

\begin{figure}[h]
 \centering
 {\includegraphics[width=0.60\textwidth]{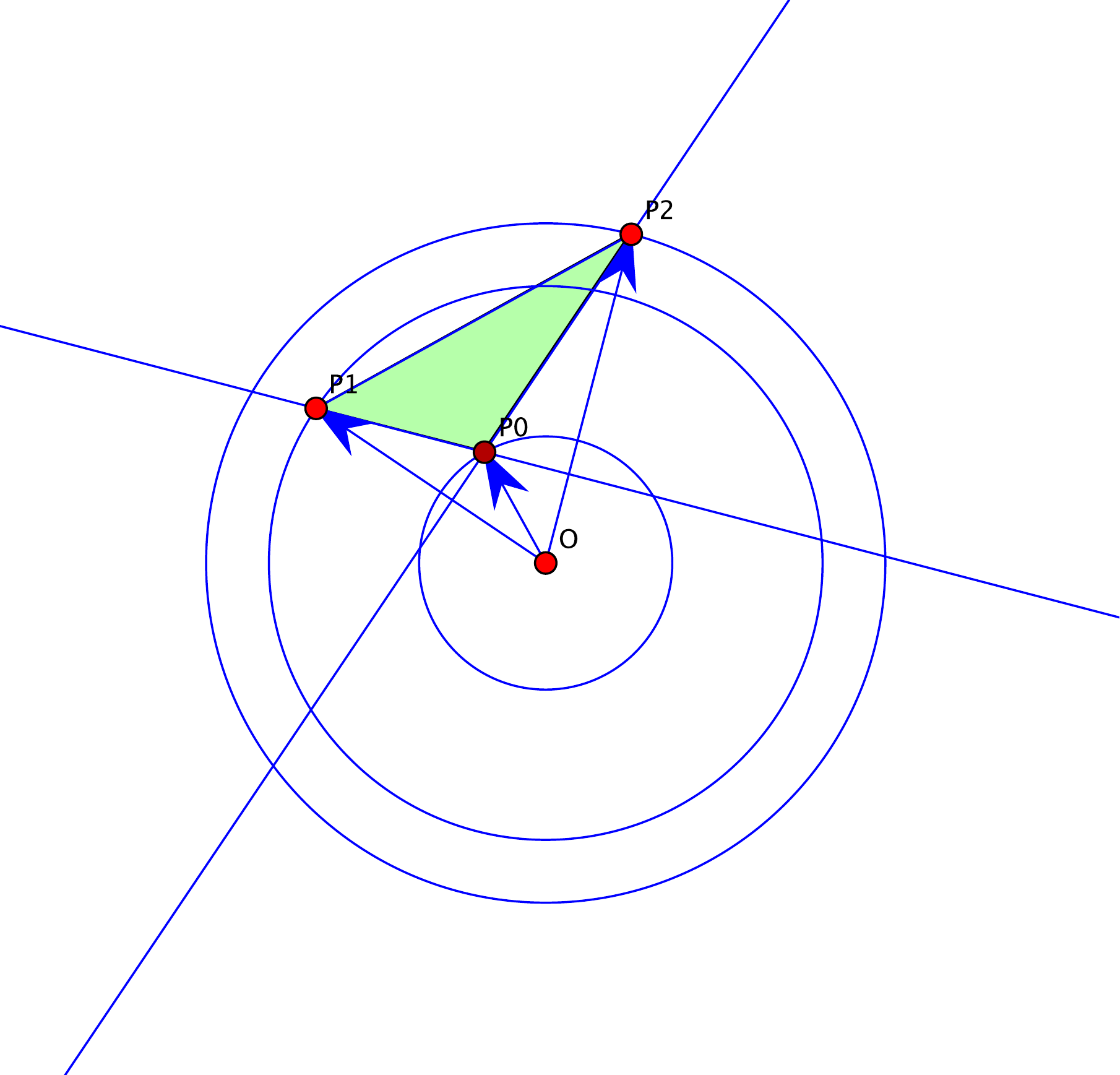}}
 \caption{ A configuration $p=(p_0,p_1,p_2)$ corresponding to a saddle point. Note that the origin $O$ is the orthocenter of triangle $p_0p_1p_2$.}
 \label{FigSaddle}
\end{figure}

\medskip \noindent
Observing that, up to $SO(d)$, a Gram matrix detrmines two reflected vector configurations
with corresponding volumes of opposite signs and equal absolute value, we may eatablish
by induction on $d\geq 2$, that the negative root corresponds to the absolute maximum and minimum, while a simple smallest positive root corresponds to reflected critical configurations which are neither local minima nor local maxima.
For $d=2$, the deformation space is a surface and a `seddle point' configuration 
is illustrated in Figure~\ref{FigSaddle}. Note that, in absolute value, a local variation of
$p_0$ increases the area of triangle $p_0p_1p_2$, while a local variatin of $p_1$ decreases
the area.

\medskip \noindent
Let us observe, for the latter part of the argument, when $s_0 < s_1\leq s_2\leq ...\leq s_d$
and $ \langle p_i,p_j\rangle = \alpha_{+}\in (s_0,s_0^{1/2}s_1^{1/2})$, $i\neq j$, that in the hyperplane passing
through $p_1,...,p_d$ we obtain a $(d-1)$-dimensional configuration with $d$ vectors $q_i$,
where

\begin{equation}\label{descent}
p_i=\mu p_0+q_i, \ \  \  \langle p_0,q_i\rangle=0, \ \ i=1,...,d
\end{equation}

\noindent
with $\mu p_0$ in the indicated hyperplane. Since $\mu=\alpha_{+}/s_0$, we find 

$$ \langle q_i,q_j\rangle  = \alpha_{+}(1-\frac{\alpha_{+}}{s_0})< 0  $$

\noindent
 which is a maximum or minimum configuration for $\langle q_i,q_i\rangle=s_i - \frac{\alpha_{+}^2}{s_0} > 0$.

\medskip \noindent
The absence of local  minima or maxima other than the absolute ones shows that 
$ \prod_{i=0}^d S^{d-1}_{\sqrt s_i}  \setminus \{ p : V(p)=0 \}   $ 
has exactly {\em two connected components} which are isomorphically exchanged by any orientation-reversing orthogonal transformation.

\begin{theorem}
Let $s_0\leq s_1\leq ,...\leq s_d$ denote the squared lengths of a complete set of edge
representatives for a $d$-dimensional diamond framework ${\cal D}(p)$.
The deformation space of this framework is a manifold of dimension ${d+1\choose 2}-1$ and can be described as the quotient $\prod_{i=0}^d S^{d-1}_{\sqrt s_i}  \setminus \{ p : V(p)=0 \} /O(d)$.
The squared volume of a fundamental parallelepiped (unit cell) has a maximum corresponding to the unique negative root $\alpha_{-}$ of (\ref{Gram}). 
\end{theorem}

\medskip
The indicated quotient space is a manifold because $O(d)$ acts without fixed points on $V(p)\neq 0$ and local charts can be obtained as transversals to orbits. An alternative argument
will follow from using Gram coordinates.

\medskip \noindent 
We have already implicated in our considerations the image of the parameter space
$\prod_{i=0}^d S^{d-1}_{\sqrt s_i}   \setminus \{ p : V(p)=0 \} /O(d)$ obtained via Gram matrices
$G(p)=(\langle p_i,p_j\rangle )_{0\leq i,j\leq d}$, with prescribed diagonal entries $\langle p_i,p_i\rangle=s_i$, $i=0,...,d$.  The image is the locus of all {\em semipositive} 
$(d+1)\times(d+1)$ symmetric matrices of rank $d$ with positive diagonal minors.

\medskip \noindent
While it is generally true  that deformation spaces of periodic frameworks are
connected components of {\em semi-algebraic sets} \cite{BS4}, the case of diamond 
frameworks allows very explicit treatment. We are going to describe another image of the
deformation space, given by the Gram marices of the generators $v_i=p_i-p_0$, $i=1,...,d$
of the corresponding periodicity lattices.

\medskip \noindent
We denote by $\omega=(\omega_{ij})_{1\leq i,j\leq d}$, the symmetric $d\times d$ matrix 
with  entries

\begin{equation}\label{omega}
\omega_{ij}=\langle v_i,v_j\rangle=\langle p_i-p_0,p_j-p_0\rangle =
 \langle p_i,p_j\rangle -\langle p_i+p_j,p_0\rangle + s_0
\end{equation}

\medskip \noindent
From $\omega_{ii}=s_i+s_0-2\langle p_i,p_0\rangle$ we retrieve all entries

\begin{equation}\label{retrieve0}
\langle p_i, p_0\rangle =  \frac{1}{2}(s_i+s_0-\omega_{ii}), \ \ i=1,...,d
\end{equation}

\noindent
and then

\begin{equation}\label{retrieve}
\langle p_i, p_j\rangle =   \omega_{ij} + \frac{1}{2}(s_i+s_j-\omega_{ii}-\omega_{jj}), 
\ \ 1\leq i \neq j \leq d
\end{equation}

\medskip \noindent
Thus, $G(p)$ and $\omega$ contain equivalent information and the {\em affine transformation} defined by (\ref{retrieve0}) and (\ref{retrieve})  take the algebraic hypersurface
$detG(p)=0$ to a degree $(d+1)$ hypersurface in the vector space of $d\times d$
 symmetric matrices (with coordinates $\omega_{ij}, 1\leq i\leq j \leq d$). Line and column operations give the equivalent expression as a `bordered determinant' of $\omega$,
namely

\begin{equation}\label{hypersurface}
det \left( \begin{array}{cccc} s_0 & \frac{1}{2}(s_1-s_0-\omega_{11}) & ... & 
\frac{1}{2}(s_d-s_0-\omega_{dd}) \\
\frac{1}{2}(s_1-s_0-\omega_{11} ) & \omega_{11} & ... & \omega_{1d} \\
.. & .. & ... & .. \\
\frac{1}{2}(s_d-s_0-\omega_{dd} ) & \omega_{d1} & ... & \omega_{dd} \end{array} \right) =0
\end{equation}

\medskip \noindent
The image of the deformation space consists of the intersection of this hypersurface with
the open cone of {\em positive definite} symmetric matrices.

\medskip \noindent
{\bf Remark:}\ Strictly speaking, we should indicate by an index $s=(s_0,...,s_d)$ the fact that the diagonal entries in the Gram matrix $G(p)=G_s(p)$ are fixed by $\langle p_i,p_i\rangle =s_i$. 
While all hypersurfaces $detG_s(p)=0$ in coordinates $a_{ij}=\langle p_i,p_j\rangle$, $0\leq i<j\leq d$ are isomorphic  to the case $s_i=1$, $i=0,...,d$ under linear transformations of the form $a_{ij}\mapsto a_{ij}/(s_is_j)^{1/2}$, what {\em varies} with $s$ is the intersection with
the locus of semipositive definite matrices with positive diagonal minors (or, in terms of
coordinates $w_{ij}$ in (\ref{hypersurface}), the intersection with the positive definite cone).

\medskip \noindent
Let us show here that for all $s$, $s_i>0, i=0,...,d$, there are no singularities of $detG_s(p)=0$
in this intersection. With $a_{ij}=\langle p_i,p_j\rangle $ as coordinates, the vanishing of the
gradient would make all off-diagonal minors equal to zero. Since $G_s(p)$ is not invertible,
the diagonal minors must vanish as well and we obtain a contradiction. The argument shows in fact that the singularity locus of $detG_s(p)=0$ is made of symmetric matrices of rank {\em
strictly less} than $d$.

\medskip \noindent
We shall investigate in more detail the case $d=2$, which involves {\em Cayley nodal cubic surfaces}.

\section{Planar diamantine frameworks}\label{planar}

For $d=2$, the deformation space of a diamantine framework is a surface. The topology of this
surface can be determined with relative ease by observing that in the quotient

$$  \prod_{i=0}^2 S^1_{\sqrt s_i} /SO(2) $$

\noindent
the action of $SO(2)=S^1$ can be used to fix $p_2$. Thus, we may use as parameters the angles
$(\phi_0,\phi_1)\in S^1\times S^1$ of $p_0$ and $p_1$ with $p_2$ and exclude from this torus,   the degenerate locus corresponding to collinear vertices $p_0,p_1,p_2$. 

\begin{figure}[h]
 \centering
 {\includegraphics[width=0.85\textwidth]{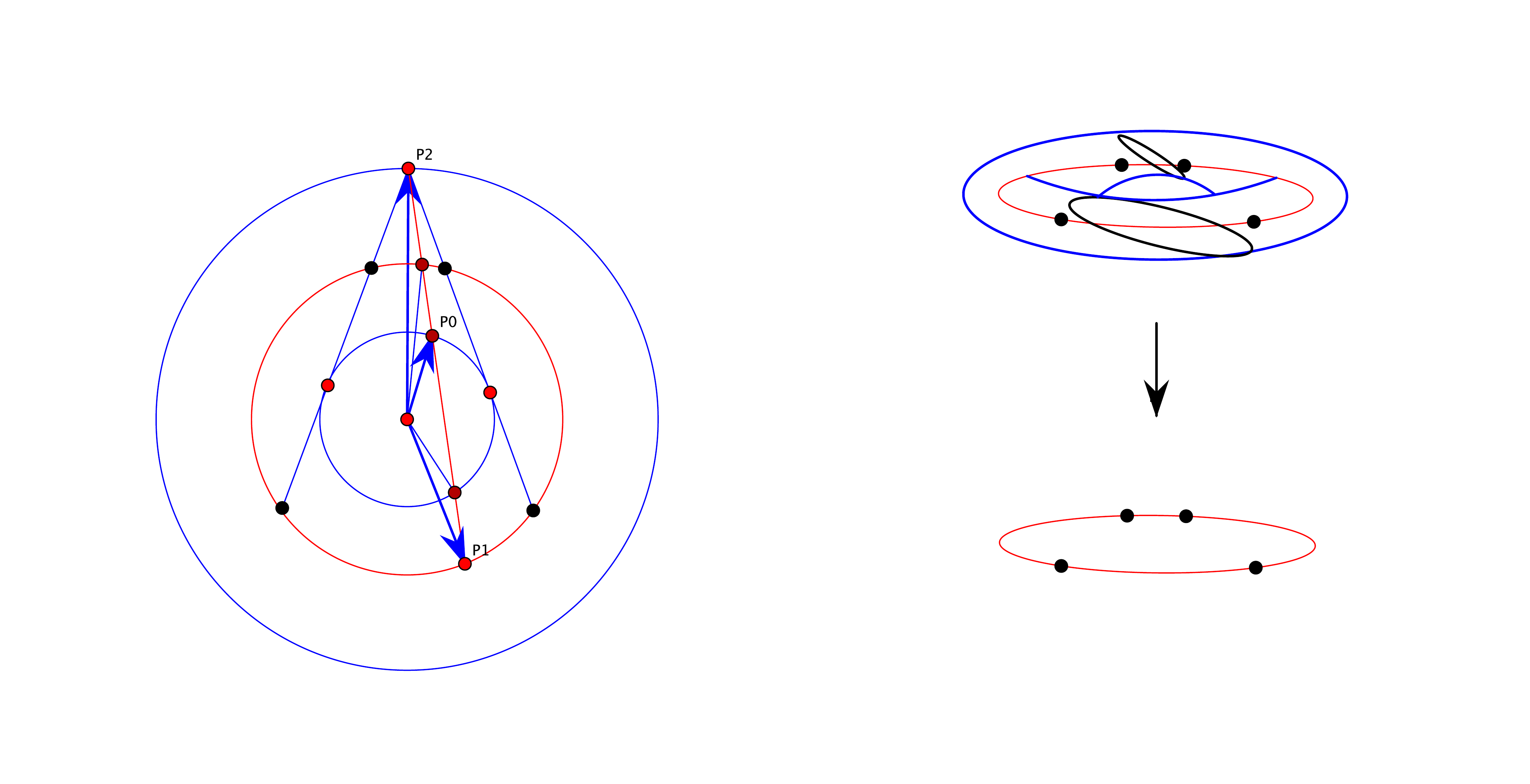}}
 \caption{ Topology determination in the generic case $s_0 < s_1 < s_2$. With $p_2$ fixed, 
the possibilities for $(p_0,p_1)$ are parametrized by the torus $S^1\times S^1$. The picture indicates the projection on the secon factor parametrizing $p_1$. The two dark  loops
depicted on the torus represent degenerate configurations with collinear vertices
 $p_0,p_1,p_2$ .}
 \label{FigGeneric}
\end{figure}

\medskip \noindent
When $s_0 < s_1\leq s_2$ the torus is cut by two disjoint loops
into two open cylinders (isomorphic by reflection in the $p_2$ axis), The generic case
is illustrated in Figure~\ref{FigGeneric}.

\medskip \noindent
For $s_0=s_1 < ls_2$
we cut along intersecting loops  which leave two (reflection isomorphic)
topological open discs in the complement. 
The case $s+0=s_1=s_2$ is the simplest, since we
exclude $\phi_0=0$, $\phi_1=0$ and $\phi_0=\phi_1$, leaving two topological open discs
in the complemeent. This gives:

\begin{prop}
The deformation space of a planar diamantine framework is an open cylinder when 
$s_0 <s_1\leq s_2$ and an open disc when $s_0=s_1\leq s_2$.
\end{prop}

\medskip \noindent
{\bf Remark:}\ The topological change from disc to cylinder is reflected in the apparition of the saddle point for the area function.

\medskip \noindent
With these preparations, we may follow the scenario described in the previous section 
whereby the deformation space is mapped to the {\em positive definite cone} (of symmetric
$2\times 2$ matrices) by associating to a framework the Gram matrix $\omega=\omega(p)$
of the two generators of the periodicity lattice $p_i-p_0$, $i=1,2$.

\begin{figure}[h]
 \centering
 {\includegraphics[width=0.60\textwidth]{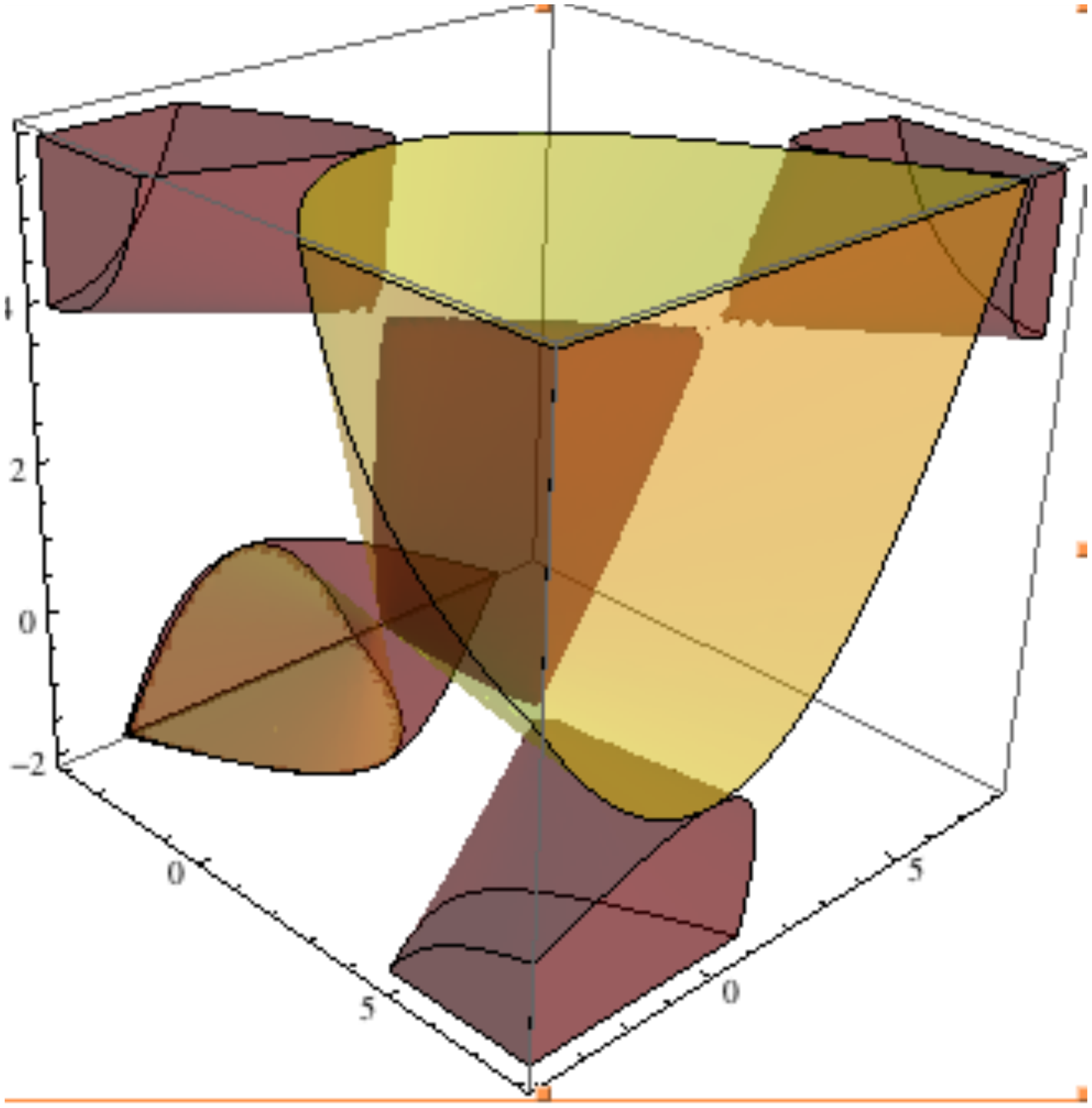}}
 \caption{ The Cayley cubic surface and the positive semidefinite cone passing through its four nodes (rendered as small gaps between nonsingular connected components) . The 2-diamond deformation surface is the part inside the cone which wraps the  tetrahedron with vertices at the nodes. The three lines of tangency along the three edges from the origin are excluded.}
 \label{FigCayleyCubic}
\end{figure}

\medskip \noindent
We illustrate the `standard' case $s_0=s_1=s_2=1$ which allows the most symmetric 
presentation, already depicted in Figure~\ref{Fig2diamond}.  Equation (\ref{hypersurface}) takes the form

$$ f(\omega)=f(\omega_{11},\omega_{12},\omega_{22})= $$

\begin{equation}\label{wCayley}
=\frac{1}{4}\omega_{11}\omega_{22}(2\omega_{12}-\omega_{11}-\omega_{22})+
\omega_{11}\omega_{22}-\omega_{12}^2=0.
\end{equation}

\medskip \noindent
The gradient vector

\begin{equation}\label{gradient}
\bigtriangledown f (\omega)=\left( \begin{array}{c}
\omega_{22}(\frac{1}{2}\omega_{12}-\frac{1}{2}\omega_{11}-\frac{1}{4}\omega_{22}+1) \\
\frac{1}{2}\omega_{11}\omega_{22}-2\omega_{12} \\
\omega_{11}(\frac{1}{2}\omega_{12}-\frac{1}{2}\omega_{22}-\frac{1}{4}\omega_{11}+1) 
\end{array} \right)
\end{equation}

\medskip \noindent
vanishes at the four points: $(0,0,0),\ ((4,0,0),\ (0,0,4),\ (4,4,4)$. At this stage, one may recognize that our cubic  (\ref{wCayley}) is projectively equivalent with the {\em Cayley nodal cubic 
surface} usually given in projective coordinates $(x_0:x_1:x_2:x_3)$ by the vanishing of the third symmetric polynomial

\begin{equation}\label{standardCayley} \sigma_3(x)=x_1x_2x_3+x_0x_2x_3+x_0x_1x_3+x_0x_1x_2 =0. 
\end{equation}

\medskip \noindent
Up to projective equivalence, , the Cayley cubic is completely characterized by the fact that it
has exactly four singularities which are all nodal points (which may be seen as the vertices of a tetrahedron whose edge-supporting lines belong to the cubic). Although equation (\ref{standardCayley}) immediately reveals projective symmetry under the permutation group
of the four coordinates $x_i$, we shall retain the description (\ref{wCayley}) in terms of $\omega$, since we have to intersect our cubic with the {\em positive definite cone} defined by

\begin{equation}\label{positiveCone}
tr(\omega )=\omega_{11}+\omega_{22} > 0 \ \ \mbox{and} \ \ 
det(\omega )=\omega_{11}\omega_{22}-\omega_{12}^2 > .0
\end{equation}

\medskip \noindent
The bounding cone $det(\omega )=0$ intersects the cubic (\ref{wCayley}) along the three double lines given by the edges  through $(0,0,0)$ of the four nodes tetrahedron. Moreover, $det(\omega)=0$ is also the tangent cone at the singularity $(0,0,0)$ for (\ref{wCayley}). 
 
\section{Auxetic trajectories}\label{aux}

The notion of {\em auxetic behaviour} is formulated in the context of elasticity theory and is
an expression of negative Poisson's ratios \cite{L, ENHR, GGLR}. Simply phrased, auxetic behaviour means becoming laterally wider when streched and thinner when compressed.
Our purely mathematical treatment of periodic framework deformations cannot implicate 
physical properties per se and we shall rely on the geometrical alternative proposed in 
\cite{BS5} and illustrated in \cite{BS4} for periodic pseudo-triangulations. 

\medskip \noindent
The auxetic property, from the greek word 
for growth, refers to certain one-parameter deformations of a periodic framework and not to the
framework itself, which may allow, in general, non-auxetic deformations as well. We use the following definition which involves the positive semidefinite cone $\bar{\Omega}(d)$ in the vector space of $d\times d$ symmetric matrices i.e. the cone made of all symmetric matrices
with non-negative eigenvalues.

\begin{definition}\label{Daux}
A smooth one-parameter deformation of a periodic framework in $R^d$ is called an auxetic
path, or simply auxetic, when the corresponding curve of Gram matrices $\omega(\tau)$ for
a set of independent generators of the periodicity lattice has all its tangent vectors $d\omega(\tau)/d\tau$ in the positive semidefinite cone of $d\times d$ symmetric matrices.
\end{definition}

\medskip \noindent
{\bf Remark.}\ The similarity with the notion of causal-line in special relativity is apparent.
A causal-line in a Minkowski space is a smooth curve with all its tangents in the {\em light cone}. For $d=2$, the boundary of the positive semidefinite cone $\bar{\Omega}(2)$ is
given by the vanishing of a single (Lorenzian) quadratic form and the two notions coincide.
For $d> 2$, the geometry of the cone $\bar{\Omega}(d)$ is more complicated \cite{Gr, B},
but the analogy persists.

\medskip
We shall make use of the description obtained in Section~\ref{planar} for the deformation surface of  2-diamnond and  investigate the possibility of auxetic trajectories on this surface.

\medskip \noindent
In order to identify the region of this surface made of points allowing some non-trivial auxetic
trajectory through them, we must look for points where some non-zero  tangent vector belongs to the positive semidefinite cone (as a free vector). The boundary of this region consists therefore
of points with tangent planes parallel to a tangent plane of the cone surface along a generating line. This locus can be detected algebraically as follows.

\medskip \noindent
We use the natural Euclidean norm $||\ ||_{tr}$  defined  on the space of $2\times 2$ symmetric matrices  by

\begin{equation}\label{trNorm}
 ||\omega||_{tr}=tr(\omega^2)^{1/2} =(\omega_{11}^2+2\omega_{12}^2+\omega_{22}^2)^{1/2}  \end{equation}

\medskip \noindent
The positive semidefinite cone is self-dual with respect to this norm \cite{Gr}, hence verifying that a tangent plane corresponts to a tangent plane of the cone $det(\omega)=0$ along a
generating line amounts to verifying that the normal direction with respect to (\ref{trNorm})
satisfies the cone equation. The required normal is immediately obtained from the gradient formula (\ref{gradient}) by halving the middle coordinate. Thus, the equation for the region's
boundary is the quartic

$$ g(\omega)=\omega_{11}\omega_{22}
(\frac{1}{2}\omega_{12}-\frac{1}{2}\omega_{11}-\frac{1}{4}\omega_{22}+1)
(\frac{1}{2}\omega_{12}-\frac{1}{2}\omega_{22}-\frac{1}{4}\omega_{11}+1)  $$

\begin{equation}\label{boundary}
-(\frac{1}{4}\omega_{11}\omega_{22}-\omega_{12})^2=0
\end{equation}

 \medskip
The curve of degree twelve resulting from the intersection of the Cayley cubic and the above quartic can be first guessed and then proven to be made of the six lines supporting the edges
of the tetrahedon of nodes, counted with multiplicity two. For instance, intersecting the quartic
(\ref{boundary}) with the plane $\omega_{12}=0$ gives the four lines resulting from the factorization 

\begin{equation}\label{factors}
g(\omega_{11},0,\omega_{22})=\frac{1}{8}\omega_{11}\omega_{22}
(\omega_{11} + \omega_{22}-4)(\omega_{11} + \omega_{22}-2)
\end{equation}

\noindent
The first three factors correspond with the three edges of the tetrahedron of nodes on the
face $\omega_{12}$. By symmetry and the fact that auxetic trajectories cannot exist in a neighborhood of the standard 2-diamond, we conclude that the region of the deformation surface consisting of points allowing a non-trivial auxetic deformation path is given by
intersecting with the open half-space

\begin{equation}\label{auxeticRegion}
 \omega_{11}-\omega_{12}+\omega_{22} <  4 
\end{equation}

\noindent
If we express (\ref{auxeticRegion}) in terms of the unit vectors $p_i$, $i=0,1,2$,  we find the
condition

\begin{equation}\label{pointed}
-1 < \langle p_0,p_1\rangle + \langle p_1,p_2\rangle + \langle p_2,p_0\rangle
\end{equation}

\begin{figure}[h]
 \centering
 {\includegraphics[width=1\textwidth]{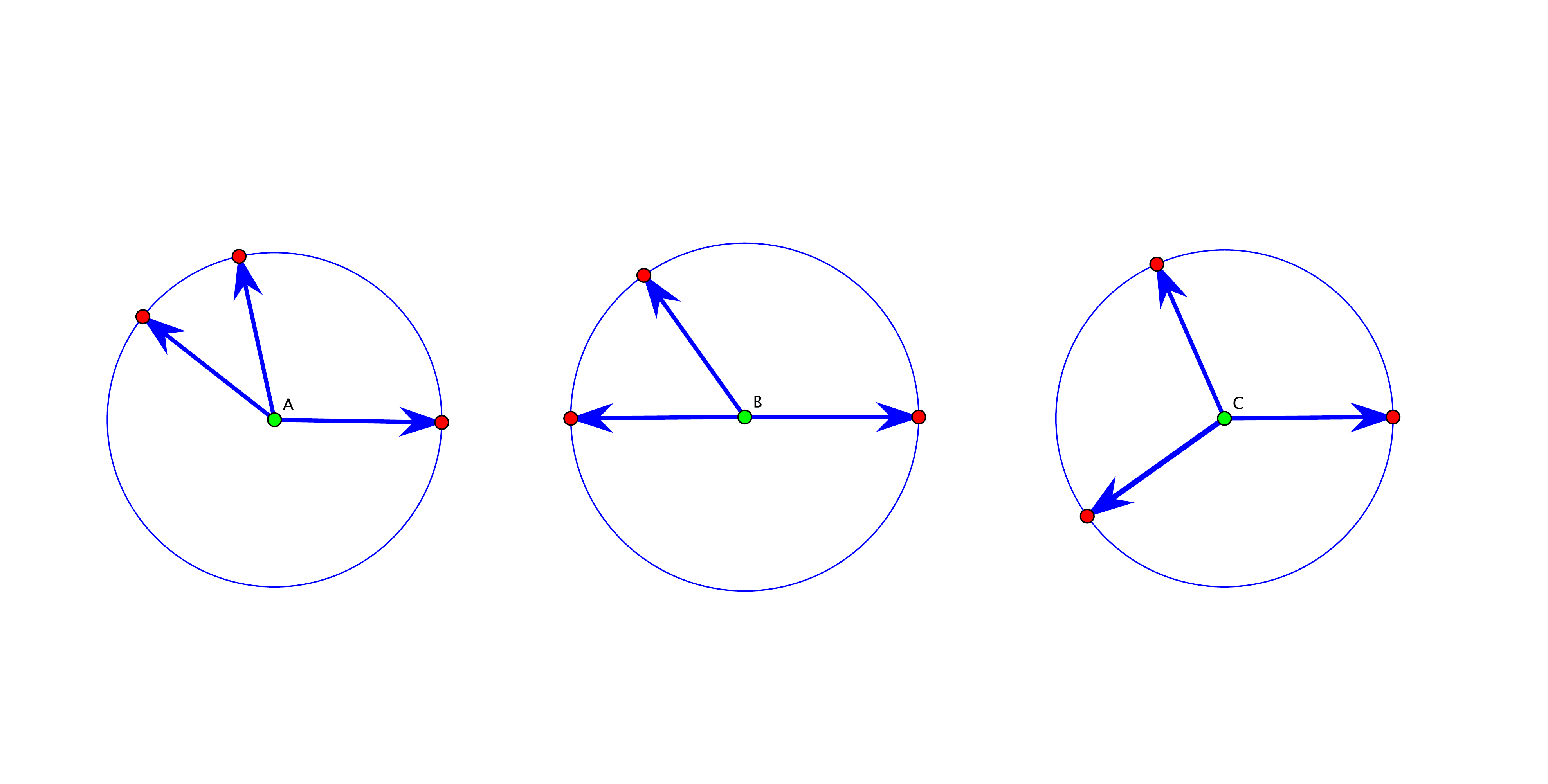}}
 \caption{ Three unit vectors in a pointed configuration (A), boundary (transition) case (B) and
non-pointed configuration (C), corresponding to $<, =$ and  $>$ in relation (\ref{pointed}). }
 \label{FigTrio}
\end{figure}

\noindent
Geometrically, this means {\em pointedness}: the three (edge) vectors $p_i$ aree all on one side
of some line through the origin.  Equivalently, one of the three sectors determined  around the origin by the three vectors has an angle larger than $\pi$ as illustrated on the left configuration
of Figure~\ref{FigTrio}.

\medskip
 We summarize our conclusion in the following statement.

\begin{prop}\label{auxeticCapability}
A framework obtained from the standard 2-diamond framework allows some non-trivial (small)
auxetic deformation path if and only if pointed at all vertices.
\end{prop}

\noindent
{\bf Remarks.}\ We have argued above for pointedness at the vertex of reference. Diamantine frameworks are vertex transitive under full crystallographic symmetry, hence pointedness at all vertices follows. The reentrant honeycomb in Figure~\ref{Figauxetic2D} is pointed at all vertices. For the relevance of pointedness in periodic pseudo-triangulations we refer to \cite{BS4,BS5,BS6}. 


\medskip

C. Borcea 

\noindent
              Department of Mathematics, Rider University, Lawrenceville, NJ 08648, USA 
              
\medskip
 I. Streinu 

\noindent
              Department of Computer Science, Smith College, Northampton, MA 01063, USA 
	
\end{document}